\documentclass[12pt]{article}
\usepackage{amsmath}
\usepackage{amssymb}
\usepackage{amsmath,amssymb,amsbsy,amsfonts,amsthm,latexsym,
amsopn,amstext,amsxtra,euscript,amscd}

\begin{document}

\def\A{\mathbb{A}}
\def\B{\mathbf{B}}
\def \C{\mathbb{C}}
\def \F{\mathbb{F}}
\def \K{\mathbb{K}}

\def \Z{\mathbb{Z}}
\def \P{\mathbb{P}}
\def \R{\mathbb{R}}
\def \Q{\mathbb{Q}}
\def \N{\mathbb{N}}
\def \Z{\mathbb{Z}}

\def\B{\mathcal B}
\def\e{\varepsilon}

\def\cA{{\mathcal A}}
\def\cB{{\mathcal B}}
\def\cC{{\mathcal C}}
\def\cD{{\mathcal D}}
\def\cE{{\mathcal E}}
\def\cF{{\mathcal F}}
\def\cG{{\mathcal G}}
\def\cH{{\mathcal H}}
\def\cI{{\mathcal I}}
\def\cJ{{\mathcal J}}
\def\cK{{\mathcal K}}
\def\cL{{\mathcal L}}
\def\cM{{\mathcal M}}
\def\cN{{\mathcal N}}
\def\cO{{\mathcal O}}
\def\cP{{\mathcal P}}
\def\cQ{{\mathcal Q}}
\def\cR{{\mathcal R}}
\def\cS{{\mathcal S}}
\def\cT{{\mathcal T}}
\def\cU{{\mathcal U}}
\def\cV{{\mathcal V}}
\def\cW{{\mathcal W}}
\def\cX{{\mathcal X}}
\def\cY{{\mathcal Y}}
\def\cZ{{\mathcal Z}}

\def\f{\frac{|\A||B|}{|G|}}
\def\AB{|\A\cap B|}
\def \Fq{\F_q}
\def \Fqn{\F_{q^n}}

\def\({\left(}
\def\){\right)}
\def\fl#1{\left\lfloor#1\right\rfloor}
\def\rf#1{\left\lceil#1\right\rceil}
\def\Res{{\mathrm{Res}}}

\newcommand{\comm}[1]{\marginpar{
\vskip-\baselineskip \raggedright\footnotesize
\itshape\hrule\smallskip#1\par\smallskip\hrule}}

\newtheorem{lem}{Lemma}
\newtheorem{lemma}[lem]{Lemma}
\newtheorem{prop}{Proposition}
\newtheorem{proposition}[prop]{Proposition }
\newtheorem{thm}{Theorem}
\newtheorem{theorem}[thm]{Theorem}
\newtheorem{cor}{Corollary}
\newtheorem{corollary}[cor]{Corollary}
\newtheorem{prob}{Problem}
\newtheorem{problem}[prob]{Problem}
\newtheorem{ques}{Question}
\newtheorem{question}[ques]{Question}
\newtheorem{rem}{Remark}

\title{Kloosterman sums in residue rings}
\author{{J.~Bourgain}\\
\normalsize{Institute for Advanced Study,}\\
\normalsize {Princeton, NJ 08540, USA}\\
\normalsize{\tt bourgain@ias.edu} \\
\and\\
{M. Z. Garaev}
\\
\normalsize{Centro de Ciencias Matem\'{a}ticas,}
\\
\normalsize{Universidad Nacional Aut\'onoma de M\'{e}xico,}\\
\normalsize{Morelia 58089, Michoac\'{a}n, M\'{e}xico}\\
\normalsize{\tt garaev@matmor.unam.mx}
 }

%\small Mathematics Subject Classifications: 11L05}

\date{\empty}

\pagenumbering{arabic}

\maketitle

\begin{abstract}

In the present paper, we generalize some of the results on Kloosterman sums proven in~\cite{BG}
for prime moduli  to general moduli. This requires to establish the corresponding additive
properties of the reciprocal set
$$
I^{-1}=\{x^{-1}:\quad x\in I\},
$$
where $I$ is an interval in the ring of residue classes modulo a large positive integer. We apply
our bounds on multilinear exponential sums to the Brun-Titchmarsh theorem and the estimate
of very short Kloosterman sums, hence generalizing our earlier work to the setting of general modulus.
\end{abstract}
\maketitle

%\newpage

\section{Introduction}

In what follows,
$\Z_m$ denotes the ring of residue classes modulo a large positive
integer $m$ which frequently will be associated with the set
$\{0,1,\ldots, m-1\}$. Given an integer $x$ coprime to $m$ (or an
invertible element of $\Z_m$) we use $x^*$ or $x^{-1}$ to denote its
multiplicative inverse modulo~$m$.

Let $I$ be an interval in $\Z_m$. In the present paper we establish some additive properties of the reciprocal-set
$$
I^{-1}=\{x^{-1}:\quad x\in I\}.
$$
We apply our results to estimate some double Kloosterman sums,
to Brun-Titchmarsh theorem and, involving  multilinear exponential sum
bounds of general modulus, we estimate short Kloosterman sums. These extends some results of our work~\cite{BG} from prime moduli
to the general.

Throughout the paper we use the abbreviation $e_m(z):=e^{2\pi i z/m}$.

\section{Statement of our results}

We start with the additive properties of reciprocal-set.
\begin{theorem}
\label{thm:kI*=kI* I=[1,N]} Let $I=[1,N]$. Then the number $J_{2k}$
of solutions of the congruence
$$
x_1^*+\ldots+x_k^*\equiv x_{k+1}^*+\ldots+x_{2k}^*\pmod m,\qquad x_1,\ldots, x_{2k}\in I,
$$
satisfies
$$
J_{2k}<(2k)^{90k^3}(\log N)^{4k^2}\Bigl(\frac{N^{2k-1}}{m}+1\Bigr)N^{k}.
$$
\end{theorem}

The following statement is a version of Theorem~\ref{thm:kI*=kI* I=[1,N]}, where
the variables $x_j$ are restricted to prime numbers. By $\cP$ we
denote the set of primes.

\begin{theorem}
\label{thm:kI*=kI* I=[1,N] prime} Let $I=[1,N]$. Then the number
$J_{2k}$ of solutions of the congruence
$$
x_1^*+\ldots+x_k^*\equiv x_{k+1}^*+\ldots+x_{2k}^*\pmod m,\qquad x_1,\ldots, x_{2k}\in I\cap\cP,
$$
satisfies
$$
J_{2k}<(2k)^{k}\Bigl(\frac{N^{2k-1}}{m}+1\Bigr)N^{k}.
$$
\end{theorem}

We recall that the incomplete Kloosterman sum is the sum of the form
$$
\sum_{x=M+1}^{M+N}e_m(ax^*+bx),
$$
where $a$ and $b$ are integers, $\gcd(a,m)=1$. Here the summation over $x$ is
restricted to $\gcd(x,m)=1$ (if the range of summation is empty, then we consider this sum
to be equal to zero). As a consequence of the
Weil bounds it is known that
$$
\Bigl|\sum_{x=1}^m e_m(ax^*+bx)\Bigr|\le \tau(m) m^{1/2},
$$
see for example~\cite[Corollary 11.12]{IwKo}. This implies that for  $N<m$
one has the bound
$$
\Bigl|\sum_{x=M+1}^{M+N}e_m(ax^*+bx)\Bigr|< m^{1/2+o(1)}.
$$
For $M=0$ and $N$ very small (that is, $N=m^{o(1)}$) these sums have
been estimated by Korolev~\cite{Kor}.

The incomplete bilinear Kloosterman sum
$$
S=\sum_{x_1=M_1+1}^{M_1+N_1}\,\sum_{x_2=M_2+1}^{M_2+N_2}\alpha_1(x_1)\alpha_2(x_2)e_m(ax_1^*x_2^*),
$$
where $\alpha_i(x_i)\in\C,\, |\alpha_i(x_i)|\le 1$, is also well
known in the literature.
When $M_1=M_2=0$ the sum $S$ (in a more general form in fact) has
been estimated by Karatsuba~\cite{Kar1, Kar2} for very short ranges
of $N_1$ and $N_2$.

Theorem~\ref{thm:kI*=kI* I=[1,N]} leads to the following improvement of
the range of applicability of Karatsuba's estimate~\cite{Kar1}.

\begin{theorem}
\label{thm:Kloost Karatsuba range} Let $I_1=[1, N_1],\,
I_2=[1,N_2]$. Then uniformly over all positive integers $k_1,k_2$
and $\gcd(a,m)=1$ we have
\begin{equation*}
\begin{split}
\Bigl|\sum_{x_1\in I_1} \sum_{x_2\in I_2}\alpha_1(x_1) &\alpha_2(x_2)e_m(ax_1^*x_2^*)\Bigr|
< (2k_1)^{\frac{45k_1^2}{k_2}}(2k_2)^{\frac{45k_2^2}{k_1}}(\log m)^{2(\frac{k_1}{k_2}+\frac{k_2}{k_1})}\times \\
\times &\Bigl(\frac{N_1^{k_1-1}}{m^{1/2}}+\frac{m^{1/2}}{N_1^{k_1}}\Bigr)^{1/(2k_1k_2)}
\Bigl(\frac{N_2^{k_2-1}}{m^{1/2}}+\frac{m^{1/2}}{N_2^{k_2}}\Bigr)^{1/(2k_1k_2)}N_1N_2.
\end{split}
\end{equation*}
\end{theorem}

Given $N_1,N_2$ we choose $k_1,k_2$ such that
$$
N_1^{2(k_1-1)}<m\le N_1^{2k_1},\qquad N_2^{2(k_2-1)}<m\le N_2^{2k_2}
$$
and the bound will be nontrivial unless both $N_1,N_2$ are within
$m^{\varepsilon}$-ratio of an element of $\{m^{\frac{1}{2l}},\l\in
\Z_{+}\}$. Thus, we have the following
\begin{cor}
\label{cor:Kloost Karatsuba rang} Let $I_1=[1, N_1],\, I_2=[1,N_2]$,
where for $i=1$ or $i=2$
$$
N_i\not\in \mathop{\bigcup}_{j\ge 1} \, [m^{\frac{1}{2j}-\varepsilon},\, m^{\frac{1}{2j}+\varepsilon}].
$$
Then
$$
\max_{(a,m)=1}\Bigl|\sum_{x_1=1}^{N_1}\sum_{x_2=1}^{N_2}\alpha_1(x_1) \alpha_2(x_2)e_m(ax_1^*x_2^*)\Bigr|<
 m^{-\delta} N_1N_2
$$
for some $\delta=\delta(\varepsilon)>0$.
\end{cor}

We shall then apply our bilinear Kloosterman sum bound to the
Brun-Titchmarsh theorem and
improve the result of Friedlander-Iwaniec~\cite{FrIw} on $\pi(x;q,a)$ as
follows:

\begin{theorem}
\label{thm:BrunTitch} Let $q~\sim x^{\theta}$, where $\theta<1$ is
close to 1. Then
$$
\pi(x;q,a)<\frac{cx}{\phi(q)\log\frac{x}{q}}
$$
with $c=2-c_1(1-\theta)^2$, for some absolute constant $c_1>0$ and
all sufficiently large $x$ in terms of $\theta$.
\end{theorem}

Recall that for $(a,q)=1$, $\pi(x;q,a)$ denotes the number of primes $p\le x$ with $p\equiv a\pmod q.$
The constant $c_1$ is effective and can be made explicit. We mention
that for primes $q$ Theorem~\ref{thm:BrunTitch} is contained in our
work~\cite{BG}.

Finally, we shall apply multilinear exponential sum bounds from~\cite{B2}
(see Lemma~\ref{lem:B2} below) to establish the following estimate of a short linear Kloosterman sums.

\begin{theorem}
\label{thm:linear sqrt log m} Let $N>m^{c}$ where $c$ is a small
fixed positive constant. Then we have the bound
$$
\max_{(a,m)=1}\Bigl|\sum_{n\le N}e_m(an^*)\Bigr|< \frac{(\log\log m)^{O(1)}}{(\log m)^{1/2}}\,N,
$$
where the implied constants may depend only on $c$.
\end{theorem}

This improves some results of Korolev~\cite{Kor}. We also refer the
reader to~\cite{Kor1} for some variants of the problem. We remark that
a stronger bound is claimed in~\cite{KarKKR}, but the proof there
is in doubt.

Since
$$
\sum_{n=1}^m e_m(an^*)=\mu(m),
$$
in Theorem~\ref{thm:linear sqrt log m} one can assume that $N<m$. We also note that the aforementioned
consequence of the Weil bounds gives a stronger estimate in the case $N>m^{1/2+c_0}$ for any fixed constant $c_0$.

\section{Lemmas}

The following result, which we state as a lemma, has been proved by
Bourgain~\cite{B2}. It is based on results from additive
combinatorics, in particular sum-product estimates. This lemma will
be used in the proof of our results  on short Kloosterman sums.
\begin{lemma}
\label{lem:B2} For all $\gamma>0$ there exist
$\varepsilon=\varepsilon(\gamma)>0,\, \tau=\tau(\gamma)>0$ and
$k=k(\gamma)\in \Z_{+}$ such that the following holds.

Let $A_1,\ldots, A_k\subset \Z_q,$ $q$ arbitrary, and assume
$|A_i|>q^{\gamma}$\, $(1\le i\le k)$ and also
$$
\max_{\xi\in \Z_{q_1}}|A_i\cap \pi_{q_1}^{-1}(\xi)|<q_1^{-\gamma}|A_i|\quad {for\,\, all} \quad q_1| q,\, q_1>q^{\varepsilon}.
$$
Then
$$
\max_{\xi\in \Z_{q}^{*}}\Bigl|\sum_{x_1\in A_1}\ldots \sum_{x_k\in A_k}e_q(\xi x_1\ldots x_k)\Bigr|<Cq^{-\tau}|A_1|\cdots |A_k|.
$$
\end{lemma}

Here, the notation $|A\cap \pi_{q_1}^{-1}(\xi)|$ can
be viewed as the number of solutions of the congruence $ x\equiv \xi
\pmod {q_1},\, x\in A.$

Clearly, the conclusion of Lemma~\ref{lem:B2} can be stated in
basically equivalent form
$$
\max_{\xi\in \Z_{q}^{*}}\sum_{x_1\in A_1}\ldots \sum_{x_{k-1}\in A_{k-1}}\Bigl|\sum_{x_k\in A_k}e_q(\xi x_1\ldots x_{k-1} x_k)\Bigr|<Cq^{-\tau}|A_1|\cdots |A_k|.
$$
Indeed applying the Cauchy-Schwarz inequality, it follows that
\begin{eqnarray*}
\begin{split}
\Bigl(\sum_{x_1\in A_1}&\ldots \sum_{x_{k-1}\in
A_{k-1}}\Bigl|\sum_{x_k\in A_k}e_q(\xi x_1\ldots x_{k-1}
x_k)\Bigr|\Bigr)^2\le \\ & |A_1|\ldots |A_{k-1}|\sum_{x_k'\in A_k}
\Bigl|\sum_{x_1\in A_1}\ldots \sum_{x_k\in A_k}e_q(\xi x_1\ldots
x_{k-1}(x_k-x_k')\Bigr|.
\end{split}
\end{eqnarray*}
We fix $x_k'\in A_k$ such that
\begin{eqnarray*}
\begin{split}
\Bigl(\sum_{x_1\in A_1}&\ldots \sum_{x_{k-1}\in
A_{k-1}}\Bigl|\sum_{x_k\in A_k}e_q(\xi x_1\ldots x_{k-1}
x_k)\Bigr|\Bigr)^2\le \\ &|A_1|\ldots |A_{k-1}||A_k|
\Bigl|\sum_{x_1\in A_1}\ldots \sum_{x_{k-1}\in A_{k-1}}\, \sum_{x_k\in A_k'}e_q(\xi x_1\ldots
x_{k-1}x_k\Bigr|,
\end{split}
\end{eqnarray*}
where $A_k'=A_k-\{x_k'\}$. Then we observe that the set $A_k'$ also
satisfies the condition of Lemma~\ref{lem:B2}.

\bigskip

We need some facts from the geometry of numbers. Recall that a
lattice in $\R^n$ is an additive subgroup of $\R^n$ generated by $n$
linearly independent vectors. Take an arbitrary convex compact and
symmetric with respect to $0$ body $D\subset\R^n$. Recall that, for
a lattice  $\Gamma\subset\R^n$ and $i=1,\ldots,n$, the $i$-th
successive minimum $\lambda_i(D,\Gamma)$ of the set $D$ with respect
to the lattice $\Gamma$ is defined as the minimal number $\lambda$
such that the set $\lambda D$ contains $i$ linearly independent
vectors of the lattice $\Gamma$. Obviously,
$\lambda_1(D,\Gamma)\le\ldots\le\lambda_n(D,\Gamma)$. We need the
following result given in~\cite[Proposition~2.1]{BHW} (see
also~\cite[Exercise~3.5.6]{TaoVu} for a simplified form that is
still enough for our purposes).

\begin{lemma}
\label{lem:latp} We have
$$
|D\cap\Gamma|\le \prod_{i=1}^n \(\frac{2i}{\lambda_i(D,\Gamma)} + 1\).
$$
\end{lemma}

Denoting, as usual, by $(2n+1)!!$ the product of all odd positive
numbers up to $2n+1$,  we get the following

\begin{cor} \label{cor:latpoints} We have
$$\prod_{i=1}^n \min\{\lambda_i(D,\Gamma),1\} \le \frac{(2n+1)!!}{
|D\cap\Gamma|}.$$
\end{cor}

\bigskip

We also need the following lemma due to Karatsuba~\cite{Kar1}.

\begin{lemma}
\label{lem:Karatsuba} The following bound holds:
\begin{equation*}
\begin{split}
\Bigl|\Bigl\{(x_1,\ldots,x_{2k})\in [1,N]^{2k}:&\quad \frac{1}{x_1}+\ldots+\frac{1}{x_k}=\frac{1}{x_{k+1}}+\ldots+\frac{1}{x_{2k}}\Bigr\}\Bigr|\\
< &(2k)^{80k^3}(\log N)^{4k^2}N^k.
\end{split}
\end{equation*}
\end{lemma}

\section{Proof of Theorems~\ref{thm:kI*=kI* I=[1,N]} and~\ref{thm:kI*=kI* I=[1,N] prime}}

First we prove Theorem~\ref{thm:kI*=kI* I=[1,N]}. It suffices to
consider the case $kN^{k-1}<m$ as otherwise the statement is
trivial. For $\lambda=0,1,\ldots,m-1$  denote
$$
J(\lambda)=\Bigl\{(x_1,\ldots,x_k)\in I^k:\quad x_1^*+\ldots+x_k^*\equiv \lambda\pmod m\Bigr\}.
$$
Let
$$
\Omega=\{\lambda\in [1, m-1]:\quad |J(\lambda)|\ge 1\}.
$$
Since $J(0)=0$, we have
$$
J_{2k}=\sum_{\lambda\in \Omega}|J(\lambda)|^2.
$$
Consider the lattice
$$
\Gamma_{\lambda}=\{(u,v)\in \Z^2:\quad \lambda u\equiv v\pmod m\}
$$
and the body
$$
D=\{(u,v)\in \R^2:\quad |u|\le N^k,\,\, |v|\le kN^{k-1}\}.
$$
Denoting by $\mu_1,\mu_2$ the consecutive minimas of the body $D$
with respect to the lattice $\Gamma_{\lambda}$, by
Corollary~\ref{cor:latpoints} it follows
$$
\prod_{i=1}^2\min\{\mu_i, 1\}\le \frac{15}{|\Gamma_{\lambda}\cap D|}.
$$
Observe that for $(x_1,\ldots,x_k)\in J(\lambda)$ one has
$$
\lambda x_1\ldots x_k\equiv x_2\ldots x_{k}+\ldots+x_1\ldots x_{k-1}\pmod m,
$$
implying
$$
(x_1\ldots x_k,\,  x_2\ldots x_{k}+\ldots+x_1\ldots x_{k-1})\in \Gamma_{\lambda}\cap D.
$$
Thus, for $\lambda\in \Omega$ we have $\mu_1\le 1$. We split the set
$\Omega$ into two subsets:
$$
\Omega'=\{\lambda\in \Omega:\quad \mu_2\le 1\},\qquad
\Omega''=\{\lambda\in \Omega:\quad \mu_2> 1\}.
$$
We have
\begin{equation}
\label{eqn:Kloost Kar Range J2k S epsilon}
\begin{split}
J_{2k} = \sum_{\lambda\in \Omega'}|J(\lambda)|^2+\sum_{\lambda\in \Omega''}|J(\lambda)|^2.
\end{split}
\end{equation}

\bigskip

{\it Case~1\/}:  $\lambda\in \Omega'$, that is $\mu_2\le 1$. Let
$(u_i,v_i)\in \mu_i D\cap\Gamma_{\lambda},\, i=1,2,$ be linearly
independent. Then
$$
0\not=u_1v_2-v_1u_2\equiv u_1\lambda u_2-u_2\lambda u_1\equiv 0\pmod m,
$$
whence
$$
\Bigl|u_1v_2-v_1u_2|\ge m.
$$
Also
$$
\Bigl|u_1v_2-v_1u_2\Bigr|\le 2k\mu_1\mu_2 N^{2k-1}\le \frac{30kN^{2k-1}}{|\Gamma_{\lambda}\cap D|}.
$$
Thus, for $\lambda\in \Omega'$,  the number $|\Gamma_{\lambda}\cap
D|$ of solutions of the congruence
$$
\lambda u\equiv v\pmod m
$$
in integers $u,v$ with $|u|\le N^k,\, |v|\le kN^{k-1}$ is bounded by
\begin{equation}
\label{eqn:Gammalambda cap D le}
|\Gamma_{\lambda}\cap D|\le \frac{30kN^{2k-1}}{m}.
\end{equation}
Note that for $\lambda\in \Omega'$ the sets
$$
\mathcal{W}_{\lambda}:=\{(u,v);\,\, (u,v)\in\Gamma_{\lambda}\cap D,\, \gcd(u,m)=1\}
$$
are pairwise disjoint.  Therefore, if we denote by
$S(u,v)$ the set of $k$-tuples $(x_1,\ldots,x_k)$ of positive
integers $x_1,\ldots,x_k\le N$ coprime to $m$ with
$$
x_1\ldots x_k=u,\quad x_2\ldots x_{k}+\ldots+x_1\ldots x_{k-1}= v,
$$
we get
$$
\sum_{\lambda\in\Omega'}|J(\lambda)|^2=\sum_{\lambda\in\Omega'}\Bigl(\sum_{\substack{(u,v)\in\Gamma_{\lambda}\cap D\\ \gcd(u,m)=1}}\,\,\sum_{(x_1,\ldots,x_k)\in S(u,v)}1\Bigr)^2.
$$
Applying the Cauchy-Schwarz
inequality and taking into account~\eqref{eqn:Gammalambda cap D le},
we get
\begin{equation}
\label{eqn:Square of J lambda over Omega'}
\sum_{\lambda\in\Omega'}|J(\lambda)|^2\le \frac{30kN^{2k-1}}{m}\sum_{\lambda\in\Omega'}\sum_{\substack{(u,v)\in\Gamma_{\lambda}\cap D\\ \gcd(u,m)=1}}\Bigl(\sum_{(x_1,\ldots,x_k)\in S(u,v)}1\Bigr)^2
\end{equation}
From the disjointness of sets $W_{\lambda}$ it follows that the
summation on the right hand side of~\eqref{eqn:Square of J lambda
over Omega'} is bounded by the number of solutions of the system of
equations
$$
\left\{\begin{array}{llll}
x_1\ldots x_k=y_1\ldots y_k,\\
x_1\ldots x_{k-1}+\ldots+x_2\ldots x_{k}=y_2\ldots y_{k}+\ldots+y_1\ldots y_{k-1},
\end{array}
\right.
$$
in positive integers $x_i,y_j\le N$ coprime to $m.$ Hence, by
Lemma~\ref{lem:Karatsuba}, it follows that
\begin{equation}
\label{eqn:sum over Omega prima}
\sum_{\lambda\in\Omega'}|J(\lambda)|^2<30k(2k)^{80k^3}(\log N)^{4k^2}\frac{N^{3k-1}}{m}.
\end{equation}

\bigskip

{\it Case~2\/}:  $\lambda\in \Omega''$, that is $\mu_2 > 1$. Then
the vectors from $\Gamma_{\lambda}\cap D$ are linearly dependent and
in particular there is some $\widehat{\lambda}\in\Q$ such that
$$
\widehat{\lambda}x_1\ldots x_k=x_2\ldots x_{k}+\ldots+x_1\ldots x_{k-1} \quad {\rm for} \quad (x_1,\ldots,x_k)\in J(\lambda).
$$
Thus,
\begin{equation*}
\begin{split}
&\sum_{\lambda\in \Omega''}|J(\lambda)|^2\le \sum_{\widehat{\lambda}\in \Q}\Bigl|\{(x_1,\ldots,x_k)\in I^k;\, \frac{1}{x_1}+\ldots+\frac{1}{x_k}=\widehat{\lambda}\Bigr|^2\\=
&\Bigl|\Bigl\{(x_1,\ldots,x_{2k})\in [1,N]^{2k}; \quad \frac{1}{x_1}+\ldots+\frac{1}{x_k}=\frac{1}{x_{k+1}}+\ldots+\frac{1}{x_{2k}}\Bigr\}\Bigr|\\
< &(2k)^{80k^3}(\log N)^{4k^2}N^k.
\end{split}
\end{equation*}
Inserting this and~\eqref{eqn:sum over Omega prima}
into~\eqref{eqn:Kloost Kar Range J2k S epsilon}, we obtain
\begin{equation*}
\begin{split}
J_{2k}<(2k)^{90k^3}(\log N)^{4k^2}\Bigl(\frac{N^{2k-1}}{m}+1\Bigr)N^{k}
\end{split}
\end{equation*}
which concludes the proof of Theorem~\ref{thm:kI*=kI* I=[1,N]}.

The proof of Theorem~\ref{thm:kI*=kI* I=[1,N] prime} follows the
same line  with the only difference that instead of
Lemma~\ref{lem:Karatsuba} one should apply the bound
\begin{equation*}
\begin{split}
\Bigl|\Bigl\{(x_1,\ldots,x_{2k})\in ([1,N]\cap \cP)^{2k};&\quad \frac{1}{x_1}+\ldots+\frac{1}{x_k}=\frac{1}{x_{k+1}}+\ldots+\frac{1}{x_{2k}}\Bigr\}\Bigr|\\
< &(2k)^{k}\Bigl(\frac{N}{\log N}\Bigr)^k.
\end{split}
\end{equation*}

\subsection{Proof of Theorem~\ref{thm:Kloost Karatsuba range}}

Let
$$
S=\sum_{x_1\in I_1}\sum_{x_2\in I_2}\alpha_1(x_1) \alpha_2(x_2)e_m(ax_1^*x_2^*).
$$
Then by H\"{o}lder's inequality
$$
|S|^{k_2}\le N_1^{k_2-1}\sum_{x_1\in I_1}\Bigl|\sum_{x_2\in I_2} \alpha_2(x_2)e_m(ax_1^*x_2^*)\Bigr|^{k_2}.
$$
Thus, for some $\sigma(x_1)\in\C,\, |\sigma(x_1)|=1$,
$$
|S|^{k_2}\le N_1^{k_2-1}\sum_{y_1,\ldots,y_{k_2}\in I_2}\Bigl|\sum_{x_1\in I_1} \sigma(x_1)e_m(ax_1^*(y_1^*+\ldots+y_{k_2}^*)\Bigr|.
$$
Again by H\"{o}lder's inequality,
$$
|S|^{k_1k_2}\le N_1^{k_1k_2-k_1}N_2^{k_1k_2-k_2}\sum_{\lambda=0}^{p-1}J_{k_2}(\lambda;N_2)\Bigl|\sum_{x_1\in I_1} \sigma(x_1)e_m(ax_1^*\lambda\Bigr|^{k_1},
$$
where $J_k(\lambda; N)$ is the number of solutions of the congruence
$$
x_1^*+\ldots+x_k^*\equiv \lambda\pmod m,\qquad x_i\in [1,N].
$$
Then applying the Cauchy-Schwarz inequality and using
$$
\sum_{\lambda=0}^{p-1}J_{k_2}(\lambda; N_2)^2=J_{2k_2}(N_2),\qquad \sum_{\lambda=0}^{p-1}\Bigl|\sum_{x_1\in I_1} \sigma(x_1)e_m(ax_1^*\lambda\Bigr|^{2k_1}\le
mJ_{2k_1}(N_1).
$$
we get
\begin{equation}
\label{eqn:th91011}
|S|^{2k_1k_2}\le pN_1^{2k_1k_2-2k_1}N_2^{2k_1k_2-2k_2}J_{2k_1}(N_1)J_{2k_2}(N_2).
\end{equation}
Applying Theorem~\ref{thm:kI*=kI* I=[1,N]}, we obtain
\begin{equation*}
\begin{split}
|S|^{2k_1k_2}\le &(2k_1)^{90k_1^3}(2k_2)^{90k_2^3}(\log N_1)^{4k_1^2}(\log N_2)^{4k_2^2}\times \\
\times &N_1^{2k_1k_2}N_2^{2k_1k_2}\Bigl(\frac{N_1^{k_1-1}}{m^{1/2}}+\frac{m^{1/2}}{N^{k_1}}\Bigr)
\Bigl(\frac{N_2^{k_2-1}}{m^{1/2}}+\frac{m^{1/2}}{N^{k_2}}\Bigr).
\end{split}
\end{equation*}
Thus,
\begin{equation*}
\begin{split}
|S|< & (2k_1)^{45k_1^2/k_2}(2k_2)^{45k_2^2/k_1}(\log m)^{2(\frac{k_1}{k_2}+\frac{k_2}{k_1})}\times \\
\times &\Bigl(\frac{N_1^{k_1-1}}{m^{1/2}}+\frac{m^{1/2}}{N^{k_1}}\Bigr)^{1/(2k_1k_2)}
\Bigl(\frac{N_2^{k_2-1}}{m^{1/2}}+\frac{m^{1/2}}{N^{k_2}}\Bigr)^{1/(2k_1k_2)}N_1N_2,
\end{split}
\end{equation*}
which finishes the proof of Theorem~\ref{thm:Kloost Karatsuba
range}.

\section{Proof of Theorem~\ref{thm:BrunTitch}}

Let  $\varepsilon$ be a positive constant very small in terms of
$\delta=1-\theta$ (say, $\varepsilon=\delta^4$). Denote
\begin{equation*}
\begin{split}
&\cA=\{n\le x; \, n\equiv a\pmod q\},\\
&\cA_d=\{n\in\cA; \, n\equiv 0\pmod d\},\\
&S(\cA,z)=|\{n\in \cA;\, (n,p)=1\,\, {\rm for}\,\, p<z, (p,q)=1\}|,\\
& r_d=|\cA_d|-\frac{x}{qd}.
\end{split}
\end{equation*}
We take $z=D^{1/2}$, where $D$
is the level of distribution. We shall define $D$ to satisfy
$$
D\sim \Bigl(\frac{x}{q}\Bigr)^{1+c\delta^2}\sim x^{\delta+c\delta^3}\sim q^{\delta+\delta^2+O(\delta^3)},
$$
where $c$ is a suitable absolute positive constant ($c=0.01$ will
do).

Take integer $k$ such that
$$
\frac{1}{2k-1}\le\frac{\delta}{2}<\frac{1}{2k-3}.
$$
%We define $D$ such that $D^{1-c\delta^2}\sim
%x^{1-\varepsilon}q^{-1},$ where $c$ is a suitable absolute positive
%constant.
Having in mind~\cite[Theorem 12.21]{FrIw1}, we consider the
factorization $D=MN$ in the form
$$
N=q^{1/(2k-1)},\quad M=\frac{D}{N}.
$$
Following the proof of~\cite[Theorem 13.1]{FrIw1} we find that
$$
S(\cA,z)\le\frac{(2+\varepsilon)x}{\phi(q)\log D}+R(M,N).
$$
Here the remainder $R(M,N)$ is estimated by
$$
R(M,N)\ll \sum_{\substack{m\le M, n\le N\\ \gcd(mn,q)=1}}\alpha_m\beta_nr_{mn},
$$
the implied constant may depend on $\varepsilon$. Our aim is to
prove the bound $R(M,N)\ll x^{1-\varepsilon}q^{-1}$. For this we may
assume that $\alpha_m,\beta_n$ are supported on dyadic intervals
$$
0.5M_1<m\le M_1,\quad 0.5N_1<n\le N_1
$$
for some $1\le M_1\le M$ and $1\le N_1\le N$ with
$M_1N_1q>x^{1-\varepsilon}.$ Then according to~\cite[p.262]{FrIw1}
we have the bound
$$
R(M,N)\ll \frac{x}{qM_1N_1}\sum_{0<|h|\le H}\sum_{m\sim M_1}\Bigl|\sum_{n\sim N_1}\gamma(h;n)e_q(ahm^*n^*)\Bigr|+\frac{x^{1-\varepsilon}}{q},
$$
where
$$
H=qM_1N_1x^{3\varepsilon-1}\le qDx^{3\varepsilon-1}\ll x^{c\delta^3+3\varepsilon}.
$$
In particular,  $\gcd(h,q)<q^{O(\delta^3)}$. Thus, for some
$\gamma(n)\in \C$ with $|\gamma(n)|\le 1$ we have
$$
R(M,N)\ll x^{3\varepsilon}\sum_{m\le M}\Bigl|\sum_{n\le N}\gamma(n)e_{q_1}(a_1m^*n^*)\Bigr|+\frac{x^{1-\varepsilon}}{q},
$$
where, say, $q^{1-\delta^2}\le q_1\le q$ and $\gcd(a_1,q_1)=1$. Then
our Theorem~\ref{thm:Kloost Karatsuba range} applied with
$k_1=k_2=k$ implies that
$$
R(M,N)\ll MN^{1-c_0/k^2}+\frac{x^{1-\varepsilon}}{q}<D^{1-c_0\delta^2}+\frac{x^{1-\varepsilon}}{q},
$$
where $c_0>0$ is an absolute constant. Therefore, from the choice
$D\sim x^{\delta+c\delta^2}$ with $0<c<0.5c_0$, we obtain
$$
S(\cA,z)<\frac{(2-c'\delta^2)x}{\phi(q)\log(x/q)}
$$
for some absolute constant $c'>0$. The result follows.

%$$
%\pi(x;q,a)<\frac{(2-c'\delta^2)x}{\phi(q)\log(x/q)}
%$$

\section{Proof of Theorem~\ref{thm:linear sqrt log m}}

The proof of Theorem~\ref{thm:linear sqrt log m} is based on
Bourgain's multilinear exponential sum bounds for general
moduli~\cite{B2}, see Lemma~\ref{lem:B2} above. We will also need a
version of Theorem~\ref{thm:Kloost Karatsuba range} on
bilinear Kloosterman sum estimates with the variables of summation
restricted to prime and almost prime numbers.

\subsection{Double Kloosterman sums with primes and almost primes}

As a consequence of Theorem~\ref{thm:kI*=kI* I=[1,N] prime} we have
the following bilinear Kloosterman sum estimate.

\begin{corollary}
\label{cor:doubleKloostWithPrimes} Let $N_1,N_2, k_1,k_2$ be
positive integers, $\gcd(a,m)=1$. Then for any coefficients
$\alpha(p),\beta(q)\in \C$ with $|\alpha(p)|, |\beta(q)|\le 1$, we
have
\begin{equation*}
\begin{split}
\Bigl|\sum_{p\le N_1}\,\, &\sum_{q\le N_2}\alpha(p) \beta(q)e_m(ap^*q^*)\Bigr|\\
< & (2k_1)^{\frac{1}{k_2}}(2k_2)^{\frac{1}{k_1}}\Bigl(\frac{N_1^{k_1-1}}{m^{1/2}}+\frac{m^{1/2}}{N_1^{k_1}}\Bigr)^{1/(2k_1k_2)}
\Bigl(\frac{N_2^{k_2-1}}{m^{1/2}}+\frac{m^{1/2}}{N_2^{k_2}}\Bigr)^{1/(2k_1k_2)}N_1N_2,
\end{split}
\end{equation*}
where the variables $p$ and $q$ of the summations are restricted to
prime numbers.
\end{corollary}

Indeed, denoting the quantity on the left hand side by $|S|$ and
following the proof of~Theorem~\ref{thm:Kloost Karatsuba range} we
arrive at the bound (see~\eqref{eqn:th91011})
$$
|S|^{2k_1k_2}\le
mN_1^{2k_1k_2-2k_1}N_2^{2k_1k_2-2k_2}J_{2k_1}(N_1)J_{2k_2}(N_2),
$$
where in our case $J_{2k}(N)$ denotes the number of solutions of the
congruence
$$
p_1^*+\ldots+p_k^*\equiv p_{k+1}^*+\ldots+p_{2k}^*\pmod m
$$
in prime numbers $p_1,\ldots,p_{2k}\le N$. The statement then
follows by the bounds for $J_{2k}(N)$ given in
Theorem~\ref{thm:kI*=kI* I=[1,N] prime}.

\begin{lemma}
\label{lem:kI*=kI* I=[1,N] almost prime equality} Let $K,L$ be large
positive integers, $2L<K$. Then uniformly over $k$ the number
$T_{2k}(K,L)$ of solutions of the diophantine equation
$$
\frac{1}{p_1q_1}+\ldots+\frac{1}{p_kq_k}=\frac{1}{p_{k+1}q_{k+1}}+\ldots+\frac{1}{p_{2k}q_{2k}}
$$
in prime numbers $p_i,q_i$ satisfying $0.5K<p_i<K$ and $q_i<L$ is
bounded by
$$
T_{2k}(K,L)<k^{4k}\Bigl(\frac{K}{\log K}\Bigr)^k\Bigl(\frac{L}{\log L}\Bigr)^{k}.
$$
\end{lemma}

The proof is straightforward. For any given $1\le i_0\le 2k$ we have
$$
\frac{p_1\ldots p_{2k}q_1\ldots q_{2k}}{p_{i_0}q_{i_0}}\equiv 0\pmod {p_{i_0}q_{i_0}}
$$
Since $p_i\not=q_j$, it follows that $p_{i_0}$ appears in the
sequence $p_1,\ldots,p_{2k}$ at least two times. Thus, the sequence
$p_1,\ldots,p_{2k}$ contains at most $k$ different prime numbers.
Correspondingly, the sequence $q_1,\ldots,q_{2k}$ contains at most
$k$ different prime numbers. Therefore, there are at most
$$
k^{2k}\Bigl(\frac{0.9K}{\log K}\Bigr)^k k^{2k}\Bigl(\frac{1.1L}{\log L}\Bigr)^k<k^{4k}\Bigl(\frac{K}{\log K}\Bigr)^k\Bigl(\frac{L}{\log L}\Bigr)^{k}
$$
possibilities for $(p_1,\ldots,p_{2k}, q_1,\ldots, q_{2k})$. The
result follows.

\smallskip

Now following the same line as the proof of
Theorems~\ref{thm:kI*=kI* I=[1,N]} and~\ref{thm:kI*=kI* I=[1,N]
prime}, with the only difference that in the course of the proof we
substitute Lemma~\ref{lem:Karatsuba} by Lemma~\ref{lem:kI*=kI*
I=[1,N] almost prime equality}, we get the following statement.

\begin{lemma}
\label{lem:kI*=kI* I=[1,N] almost prime congruence} Let $K,L$ be
large positive integers, $2L<K$. Then uniformly over $k$ the number
$J_{2k}(K,L)$ of solutions of the diophantine equation
$$
\frac{1}{p_1q_1}+\ldots+\frac{1}{p_kq_k}\equiv \frac{1}{p_{k+1}q_{k+1}}+\ldots+\frac{1}{p_{2k}q_{2k}}\pmod {m}
$$
in prime numbers $p_i,q_i$ satisfying $0.5K<p_i<K$ and $q_i<L$ is
bounded by
$$
J_{2k}(K,L)<k^{4k}\Bigl(\frac{(KL)^{2k-1}}{m}+1\Bigr)(KL)^{k}.
$$
\end{lemma}

From Lemma~\ref{lem:kI*=kI* I=[1,N] almost prime congruence} we get
the following corollary.

\begin{corollary}
\label{cor:doubleKloostWithAlmostPrimes} Let $N,K,L, k_1,k_2$ be
positive integers, $2L<K$. Then for any coefficients
$\alpha(p),\beta(q;r)\in \C$ with $|\alpha(p)|, |\beta(q;r)|\le 1$,
we have
\begin{equation*}
\begin{split}
&\max_{\gcd(a,m)=1}\Bigl|\sum_{p\le N}\,\,\sum_{0.5K<q\le K}\,\sum_{r\le L}\alpha(p) \beta(q;r)e_m(ap^*q^*r^*)\Bigr|\\
&<  k_1^{\frac{2}{k_2}}k_2^{\frac{2}{k_1}}\Bigl(\frac{N^{k_1-1}}{m^{1/2}}+\frac{m^{1/2}}{N^{k_1}}\Bigr)^{1/(2k_1k_2)}
\Bigl(\frac{(KL)^{k_2-1}}{m^{1/2}}+\frac{m^{1/2}}{(KL)^{k_2}}\Bigr)^{1/(2k_1k_2)}NKL,
\end{split}
\end{equation*}
where the variables $p, q$ and $r$ of the summations are restricted
to prime numbers.
\end{corollary}

\subsection{Proof of Theorem~\ref{thm:linear sqrt log m}}

Denote $\varepsilon:=\log N/\log m>c$. As we have mentioned before, we can assume that $\varepsilon<4/7$.

In what follows, $r$ is a large
absolute integer constant. More explicitly, we define $r$ to be the
choice of $k$ in Lemma~\ref{lem:B2} with, say, $\gamma=1/10$.
Denote
$$
\cG=\{x<N:\quad p_1\ge N^{\alpha},\, p_r\ge N^{\beta},\, p_1p_2\ldots p_r<N^{1-\beta}\},
$$
where $p_1\ge p_2\ge\ldots\ge  p_r$ are the largest prime factors of
$x$ and
$$
0.1>\alpha> \beta>\frac{1}{\log N}
$$
are parameters to specify. Note that the number of positive integers
not exceeding $N$ and consisting on products of at most $r-1$ prime
numbers is estimated by
\begin{eqnarray*}
\sum_{k=1}^{r-1}\sum_{\substack{p_1\ldots p_k\le N\\ p_1\ge\ldots\ge p_k}} 1\ll
\frac{N}{\log N}+\sum_{k=2}^{r-1}\sum_{p_2\ldots p_k\le N^{(k-1)/k}}\frac{N}{p_2\ldots p_k\log (N/(p_2\ldots p_k))}\\
\ll \frac{N}{\log N}+\sum_{k=2}^{r-1}\sum_{p_2\le N}\ldots\sum_{p_k\le N}\frac{N}{p_2\ldots p_k\log N}\\
\ll \frac{N(\log\log N)^{r-1}}{\log N}.
\end{eqnarray*}
Here and below the implied constants may depend only on $r$. Hence,
we have
\begin{equation*}
\begin{split}
N-|\cG| \le  \frac{cN(\log\log N)^{r-1}}{\log N}+\sum_{\substack{x<N\\ p_1<N^{\alpha}}}1+\sum_{\substack{x<N\\ p_r<N^{\beta}}}1+
\sum_{\substack{x<N\\ p_1p_2\ldots p_{r}>N^{1-\beta}}}1,
\end{split}
\end{equation*}
for some constant $c=c(r)>0$. Next, we have
\begin{equation*}
\begin{split}
\sum_{\substack{x<N\\ p_1p_2\ldots p_{r}>N^{1-\beta}}}1&\le \sum_{\substack{y<N^{\beta}\\ p_1p_2\ldots p_{r}<N/y}}1\\&\ll
\sum_{\substack{y<N^{\beta}}}\, \sum_{p_2\ldots p_r<(N/y)^{(r-1)/r}}\frac{N}{y p_2\ldots p_r\log (N/(y p_2\ldots p_r))}\\
&\ll \sum_{y<N^{\beta}}\frac{N(\log\log N)^{r-1}}{y\log N}\ll \beta N (\log\log N)^{r-1}.
\end{split}
\end{equation*}
Let $\Psi(x,y)$, as usual, denote the number of positive integers
$\le x$ having no prime divisors $>y$. Thus, we have
\begin{equation*}
\begin{split}
N-|\cG| \le c_1\beta N (\log\log N)^{r-1}+\Psi(N,N^{\alpha}) +\sum_{\substack{x<N\\ p_r<N^{\beta}}}1,
\end{split}
\end{equation*}
for some constant $c_1=c_1(r)>0$.

Letting $0.1>\beta_1>\beta$ be another parameter, we similarly
observe that
\begin{equation*}
\begin{split}
\sum_{\substack{x<N\\ p_1\ldots p_{r-1}>N^{1-\beta_1}}}1 &\le \sum_{y<N^{\beta_1}}\sum_{p_1\ldots p_{r-1}\le N/y}1\\
&\ll\sum_{y<N^{\beta_1}}\frac{N(\log\log N)^{r-2}}{y\log N}
\ll \beta_1 N (\log\log N)^{r-2}.
\end{split}
\end{equation*}
Hence,
\begin{equation*}
\begin{split}
N-|\cG| \le c_1\beta N (\log\log N)^{r-1}+c_2\beta_1 N (\log\log N)^{r-2}\\ +\Psi(N,N^{\alpha})
+\sum_{\substack{x<N\\ p_r<N^{\beta}\\ p_1\ldots p_{r-1}\le N^{1-\beta_1}}}1,
\end{split}
\end{equation*}
Observing that
$$
\sum_{\substack{x<N\\ p_r<N^{\beta}\\ p_1\ldots p_{r-1}\le N^{1-\beta_1}}}1\le
\sum_{p_1\ldots p_{r-1}\le N^{1-\beta_1}}\Psi\Bigl(\frac{N}{p_1\ldots p_{r-1}}, N^{\beta}\Bigr),
$$
we get
\begin{equation*}
\begin{split}
N-|\cG| \le c_1\beta N (\log\log N)^{r-1}+c_2\beta_1 N (\log\log N)^{r-2}+\\ \Psi(N,N^{\alpha})
+\sum_{p_1\ldots p_{r-1}\le N^{1-\beta_1}} \Psi\Bigl(\frac{N}{p_1\ldots p_{r-1}}, N^{\beta}\Bigr).
\end{split}
\end{equation*}
 By the classical result of de Bruijn~\cite{Bruin} if
$y>(\log x)^{1+\delta}$, where $\delta>0$ is a fixed constant, then
$$
\Psi(x,y)\le xu^{-u(1+o(1))} \quad {\rm as}\quad u=\frac{\log x}{\log y}\to\infty.
$$
We now take
$$
\alpha=\frac{1}{\log\log m},\qquad \beta= \frac{\log\log m}{(\log m)^{1/2}},\qquad \beta_1=\beta \log\log m=\frac{(\log\log m)^{2}}{(\log m)^{1/2}}
$$
and then have
\begin{equation*}
\begin{split}
N-|\cG| & <\alpha^{\frac{1}{2\alpha}}N+\sum_{p_1\ldots p_{r-1}<N^{1-\beta_1}}\frac{N}{p_1\ldots p_{r-1}}\Bigl(\frac{\beta}{\beta_1}\Bigr)^{\frac{\beta_1}{2\beta}}+
c\beta N(\log\log m)^{r-1}\\
&<\Bigl(\alpha^{\frac{1}{2\alpha}}+(\log\log N)^{r-1}\Bigl(\frac{\beta}{\beta_1}\Bigr)^{\frac{\beta_1}{2\beta}}+c_3\beta (\log\log m)^{r-1}\Bigr)N\\
&<c_4\beta(\log\log m)^{r-1} N.
\end{split}
\end{equation*}
Therefore
\begin{equation}
\label{eqn:sum short x<N}
\Bigl|\sum_{x<N}e_m(ax^*)\Bigr|\le c_4\beta(\log\log m)^{r-1} N+\Bigl|\sum_{x\in \cG}e_m(ax^*)\Bigr|.
\end{equation}
The sum $\sum\limits_{x\in \cG}e_m(ax^*)$ may be bounded by
\begin{equation}
\label{eqn:cuadruple sum}
{\sum_{p_1}}{\sum_{p_2}}\ldots{\sum_{p_r}}\Bigl|\sum_{y}e_m(ap_1^*p_2^*\ldots p_r^*y^*)\Bigr|,
\end{equation}
where the summations are taken over primes $p_1, p_2, \ldots, p_r$
and integers $y$ such that
\begin{equation}
\label{eqn:primes p1p2p3}
p_1\ge p_2\ge \ldots \ge p_r;\quad p_1\ge N^{\alpha}; \quad p_r\ge N^{\beta};\quad p_1p_2\ldots p_r\le N^{1-\beta}
\end{equation}
and
$$
y<\frac{N}{p_1p_2\ldots p_r};\qquad P(y)\le p_r.
$$
Note that if $t$ and $T$ are such that
\begin{equation}
\label{eqn:t and T}
\Bigl(1-\frac{c_0}{\log m}\Bigr)p_r<t< p_r,\qquad  \Bigl(1-\frac{c_0}{\log m}\Bigr)\frac{N}{p_1p_2\ldots p_r}<T<\frac{N}{p_1p_2\ldots p_r},
\end{equation}
where $c_0>0$ is any constant, then we can substitute the condition on
$y$ with
\begin{equation}
\label{eqn:yyy}
P(y)\le t; \qquad y<T
\end{equation}
by changing the sum~\eqref{eqn:cuadruple sum} with an additional
term of size at most
$$
\frac{N(\log\log m)^{O(1)}}{\log m}.
$$
%Thus, for any $t$ and $T$ satisfying~\eqref{eqn:t and T} we have
%$$
%\Bigl|\sum_{x\in \cG}e_m(ax^*)\Bigr|<\frac{N(\log\log m)^{O(1)}}{\log m}+
%{\sum_{p_1}}{\sum_{p_2}}\ldots {\sum_{p_r}}\Bigl|\sum_{y}e_m(ap_1^*p_2^*\ldots p_r^*y^*)\Bigr|,
%$$
%where the summations are taken over primes $p_1, p_2, \ldots, p_r$
%and integers $y$ satisfying~\eqref{eqn:primes p1p2p3}
%and~\eqref{eqn:yyy}.

Now we split the range of summation of primes $p_1,p_2, \ldots,
p_r$ into subintervals of the form $[L, L+L(\log m)^{-1}]$ and
choosing suitable $t$ and $T$ we obtain that for some numbers
$M_1,M_2,\ldots, M_r$ with
\begin{equation}
\label{eqn:M1M2M3 Kloost restrictions}
M_1>M_2>\ldots >M_r,\quad M_1\ge \frac{N^{\alpha}}{2}, \quad M_r\ge \frac{N^{\beta}}{2},\quad M_1M_2\ldots M_r<N^{1-\beta}
\end{equation}
one has
\begin{equation}
\label{eqn:sum x in G Kloost}
\begin{split}
\Bigl|\sum_{x\in \cG}e_m(ax^*)\Bigr|&<\frac{N(\log\log m)^{O(1)}}{\log m}\\&+
(\log m)^{3r}{\sum_{p_1\in I_1}}{\sum_{p_2\in I_2}}\ldots {\sum_{p_r\in I_r}}\Bigl|\sum_{\substack{y\le M\\ P(y)\le M_r}}e_m(ap_1^*p_2^*\ldots p_r^*y^*)\Bigr|,
\end{split}
\end{equation}
where
$$
I_j=\Bigl[M_j, M_j+\frac{M_j}{\log m}\Bigr], \qquad M=\frac{N}{M_1M_2\ldots M_r}\ge N^{\beta}.
$$
Denote
$$
W=\sum_{p_1\in I_1}\sum_{p_2\in I_2}\ldots \sum_{p_r\in I_r}\Bigl|\sum_{\substack{y\le M\\ P(y)\le M_r}}e_m(ap_1^*p_2^*\ldots p_r^*y^*)\Bigr|.
$$
Applying the Cauchy-Schwarz inequality, we get
$$
W^2\le M_1M_2\ldots M_r\sum_{y\le M}\sum_{z\le M}\Bigl|\sum_{p_1\in I_1}\sum_{p_2\in I_2}\ldots \sum_{p_r\in I_r}e_m\Bigl(ap_1^*p_2^*\ldots p_r^*(y^*-z^*)\Bigr)\Bigr|.
$$
Taking into account the contribution from the pairs $y$ and $z$
with, say,
$$
\gcd(y-z,m)>e^{10\log m/\log\log m}
$$
and then fixing the pairs  $y$ and $z$ with $ \gcd(y-z,m)\le
e^{10\log m/\log\log m}$, we get the bound
\begin{equation}
\label{eqn:W^2 Kloost}
W^2\le \frac{N^2}{M}+\frac{N^2}{e^{\log m/\log\log m}}+NM|S|\le 2N^{2-\beta}+\frac{N^2}{M_1M_2\ldots M_r}|S|,
\end{equation}
where
$$
|S|=\Bigl|\sum_{p_1\in I_1}\sum_{p_2\in I_2}\ldots \sum_{p_r\in I_r}e_{m_1}(bp_1^*p_2^*\ldots p_r^*)\Bigr|.
$$
Here $b$ and $m_1$ are some positive integers satisfying
$$
\gcd(b,m_1)=1,\quad m_1\ge me^{-10\log m/\log\log m}.
$$

We consider two cases, depending on whether $M_r> N^{\alpha^3}$ or
$M_r\le N^{\alpha^3}$.

\smallskip

{\it Case 1}. Let $M_r > N^{\alpha^3}$. Hence $M_j>N^{\alpha^3}$ for
all $j=1,2,\ldots,r.$ The idea is to use Theorem~\ref{thm:kI*=kI*
I=[1,N] prime} and amplify each of these factors  to size
$m^{1/3+o(1)}$ say and then apply Lemma~\ref{lem:B2}.

Let $k_1,\ldots,k_r$ be positive integers defined from
$$
M_i^{2k_i-1}<m_1\le M_i^{2k_i+1}.
$$
Since
$M_i>N^{\alpha^3}>m^{c\alpha^3}$, it follows that
$$
k_i<\frac{1}{c\alpha^3}=\frac{(\log\log m)^3}{c}.
$$

Consequently applying H\"older's inequality, we get the bound
$$
|S|^{2^rk_1k_2\ldots k_{r}}\le \Bigl(\prod_{i=1}^{r}M_i^{2^rk_1\ldots k_r-2k_i}\Bigr)\sum_{\substack{p_{11},\ldots, p_{1k_1}\in I_1\cap\cP\\ q_{11},\ldots, q_{1k_1}\in I_1\cap\cP}}\ldots \sum_{\substack{p_{r1},\ldots, p_{r k_r}\in I_r\cap\cP\\ q_{r1},\ldots,q_{rk_r}\in I_r\cap\cP}}e^{2\pi i b\{...\}/m_1},
$$
where  $\{...\}$ indicates the expression
$$
(p_{11}^*+\ldots+p_{1k_1}^*-q_{11}^*-\ldots-q_{1k_1}^*)\cdots (p_{r1}^*+\ldots+p_{rk_r}^*-q_{r1}^*-\ldots-q_{rk_r}^*)
$$
Next, we can fix the variables $q_{ij}$ and then get that for some
integers $\mu_1,\ldots,\mu_r$ there is the bound
\begin{equation}
\label{eqn:Case 1 S via S1}
\frac{|S|}{ M_1M_2\ldots M_r}\le\Bigl(\frac{|S_1|}{M_1^{k_1}M_2^{k_2}\ldots M_r^{k_r}}\Bigr)^{1/(2^rk_1 k_2 \ldots k_r)},
\end{equation}
where
$$
S_1=\sum_{p_{11},\ldots, p_{1k_1}\in I_1\cap\cP}\ldots \sum_{p_{r1},\ldots, p_{r k_r}\in I_r\cap\cP}e^{2\pi i b(p_{11}^*+\ldots+p_{1k_1}^*-\mu_1)\cdots (p_{r1}^*+\ldots+p_{rk_r}^*-\mu_r)/m_1}.
$$
Let $A_1,\ldots, A_r$ be subsets of $\Z_{m_1}$ defined by
\begin{eqnarray*}
A_1=\{p_{11}^*+\ldots+p_{1k_1}^*-\mu_1; \quad (p_{11},\ldots,p_{1k_1})\in (I_1\cap\cP)^{k_1}\},\\
\centerline{\ldots}\\
A_r=\{p_{r1}^*+\ldots+p_{rk_r}^*-\mu_r; \quad (p_{r1},\ldots,p_{rk_r})\in (I_r\cap\cP)^{k_r}\},
\end{eqnarray*}
where $p_{ij}^*$ are calculated modulo $m_1$. Then we have
$$
S_1=\sum_{\lambda_1\in A_1}\ldots \sum_{\lambda_r\in A_r} I_1(\lambda_1)\ldots I_r(\lambda_r) e^{2\pi i b \lambda_1\ldots \lambda_r/m_1},
$$
where $I_j(\lambda)$ is the number of solutions of the congruence
$$
p_{j1}^*+\ldots+p_{jk_j}^*-\mu_j\equiv\lambda\pmod {m_1}; \qquad (p_{j1},\ldots,p_{jk_j})\in (I_j\cap\cP)^{k_j}.
$$
We apply Cauchy-Schwarz inequality to the sum over
$\lambda_1,\ldots,\lambda_{r-1}$ and get
$$
|S_1|^{2}\le J_{2k_1}(M_1)\ldots J_{2k_r}(M_r)\sum_{\lambda_1\in A_1}\ldots \sum_{\lambda_{r-1}\in A_{r-1}}\Bigl|\sum_{\lambda_r\in A_r} I_r(\lambda_r) e^{2\pi i b \lambda_1\ldots \lambda_{r-1} \lambda_r/m_1}\Bigr|^2,
$$
where
$$
J_{2k_j}(M_j)=\sum_{\lambda\in A_j} (I_j(\lambda))^2.
$$
Changing the order of summation, we get
\begin{eqnarray*}
\begin{split}
|S_1|^{2}\le & J_{2k_1}(M_1)\ldots J_{2k_{r-1}}(M_{r-1})\times \\ & \sum_{\lambda_r,\lambda_r'\in A_r} I_r(\lambda_r)I_r(\lambda_r')\Bigl|\sum_{\lambda_1\in A_1}\ldots \sum_{\lambda_{r-1}\in A_{r-1}} e^{2\pi i b\lambda_1\ldots \lambda_{r-1} (\lambda_r-\lambda_r')/m_1}\Bigr|.
\end{split}
\end{eqnarray*}
We apply the Cauchy-Schwarz inequality to the sum over
$\lambda_r,\lambda_r'$ and get
\begin{eqnarray*}
\begin{split}
|S_1|^{4}\le &(J_{2k_1}(M_1)\ldots J_{2k_{r}}(M_{r}))^2\times \\ & \sum_{\lambda_r,\lambda_r'\in A_r} \Bigl|\sum_{\lambda_1\in A_1}\ldots \sum_{\lambda_{r-1}\in A_{r-1}} e^{2\pi i b \lambda_1\ldots \lambda_{r-1} (\lambda_r-\lambda_r')/m_1}\Bigr|^2.
\end{split}
\end{eqnarray*}
We can fix $\lambda_r'\in A_r$ such that
\begin{eqnarray*}
\begin{split}
|S_1|^{4}\le & (J_{2k_1}(M_1)\ldots J_{2k_{r}}(M_{r}))^2 |A_r|\times \\ & \sum_{\lambda_r\in A_r'} \Bigl|\sum_{\lambda_1\in A_1}\ldots
\sum_{\lambda_{r-1}\in A_{r-1}} e^{2\pi i b \lambda_1\ldots \lambda_{r-1} \lambda_r/m_1}\Bigr|^2,
\end{split}
\end{eqnarray*}
where $A_r'=A_r-\{\lambda_r'\}$. Using the trivial bound
\begin{eqnarray*}
\begin{split}
\Bigl|\sum_{\lambda_1\in A_1}&\ldots
\sum_{\lambda_{r-1}\in A_{r-1}} e^{2\pi i b \lambda_1\ldots \lambda_{r-1} \lambda_r/m_1}\Bigr|^2 \\ & \le |A_1|\ldots|A_{r-1}| \Bigl|\sum_{\lambda_1\in A_1}\ldots
\sum_{\lambda_{r-1}\in A_{r-1}} e^{2\pi i b \lambda_1\ldots \lambda_{r-1} \lambda_r/m_1}\Bigr|,
\end{split}
\end{eqnarray*}
we get
\begin{eqnarray*}
\begin{split}
|S_1|^{4}\le (J_{2k_1}(M_1)& \ldots  J_{2k_{r}}(M_{r}))^2 |A_1|\ldots |A_r|\times \\ & \sum_{\lambda_r\in A_r'} \Bigl|\sum_{\lambda_1\in A_1}\ldots
\sum_{\lambda_{r-1}\in A_{r-1}} e^{2\pi i b \lambda_1\ldots \lambda_{r-1} \lambda_r/m_1}\Bigr|,
\end{split}
\end{eqnarray*}
From the definition of $A_i$ we have $|A_i|\le M_i^{k_i}$ . From the
choice of $k_i$ and Theorem~\ref{thm:kI*=kI* I=[1,N] prime} we also
have
$$
J_{2k_i}(M_i)<2(2k_i)^{k_i}M_i^{k_i}.
$$
Thus,
\begin{eqnarray}
\label{eqn:S1^4<}
|S_1|^{4}\le \Bigl(\prod_{i=1}^r(4k_i)^{k_i}M_i^{3k_i}\Bigr)\times \sum_{\lambda_r\in A_r'} \Bigl|\sum_{\lambda_1\in A_1}\ldots
\sum_{\lambda_{r-1}\in A_{r-1}} e^{2\pi i b \lambda_1\ldots \lambda_{r-1} \lambda_r)/m_1}\Bigr|,
\end{eqnarray}

Let $\gamma=1/10$ and define
$\varepsilon=\varepsilon(\gamma)>0$ to be the absolute constant from
Lemma~\ref{lem:B2}. We shall verify that the sets $A_1,\ldots,A_r$
satisfy the condition of Lemma~\ref{lem:B2} with $q=m_1$ (note that
if $A_r$ satisfies the condition of Lemma~\ref{lem:B2} then also
does $A_r'$).  From the definition of $A_i$ and the connection
between the cardinality of a set and the corresponding additive
energies, we have
\begin{equation}
\label{eqn:Ai greater than Mi to the ki}
|A_i|\ge \frac{(M_i/(2\log M_i))^{2k_i}}{J_{2k_i}(M_i)}\ge \frac{M_i^{k_i}}{2(2k_i)^{k_i}(2\log M_i)^{2k_i}}.
\end{equation}
From the choice of $k_i$ it then follows that
$$
|A_i|\ge \frac{m_1^{1/3}}{2(2k_i)^{k_i}(2\log M_i)^{2k_i}}=m_1^{1/3+o(1)}.
$$

Thus, the first condition $|A_i|>m_1^{1/10}$ is satisfied.

Next, let $q_1|m_1, q_1>m_1^{\varepsilon}$ and let $\xi\in\Z_{q_1}$.
Let $T_i$ be the number of solutions of the congruence
$$
x\equiv \xi\pmod {q_1};\qquad x\in A_i.
$$
It follows that $T_i$ is bounded by the number of solutions of the
congruence
$$
p_{1}^*+\ldots+p_{k_i}^*\equiv \xi+\mu_1 \pmod {q_1}; \quad (p_{1},\ldots,p_{k_i})\in (I_i\cap\cP)^{k_i}.
$$
Consider two possibilities here. If $M_i>q_1^{1/8}$ say, then we fix
$p_2,\ldots,p_{k_i}$ and we have at most $1+M_iq_1^{-1}$ possibilities for $p_1$. Thus, using~\eqref{eqn:Ai greater than Mi to the ki},
we get
$$
\Bigl(1+\frac{M_i}{q_1}\Bigr)M_i^{k_i-1}<\frac{M_i^{k_i}}{q_1^{1/9}}<q_1^{-1/10}|A_i|.
$$
Therefore, in this case the condition of Lemma~\ref{lem:B2} is
satisfied.

Let now $M_i<q_1^{1/8}$. Define $k_i'$ from the
condition
$$
M_i^{4k_i'+1}<q_1< M_i^{4k_i'+5}.
$$
We then have $2k_i'<k_i$. Thus,
$$
T_i\le M_i^{k_i-2k_i'}J_{2k_i'}(M_i),
$$
where $J_{2k_i'}(M_i)$, as before, denotes the number of solutions
of the congruence
$$
p_{1}^*+\ldots+p_{k_i'}^*\equiv p_{k_i'+1}^*+\ldots+p_{2k_i'}^* \pmod {q_1}; \quad (p_{1},\ldots,p_{2k_i'})\in (I_i\cap\cP)^{2k_i'}.
$$
From the choice of $k_i$ and Theorem~\ref{thm:kI*=kI* I=[1,N] prime}
we get that
$$
J_{2k_i'}(M_i)<2(2k_i)^{k_i}M_i^{k_i'}.
$$
Therefore, using~\eqref{eqn:Ai greater than Mi to the ki}
$$
T_i\le 2(2k_i)^{k_i} M_i^{k_i-k_i'}\le 2(2k_i)^{k_i} M_i^{k_i}q^{-1/9}<q_1^{-1/10}|A_i|.
$$

Thus, the condition of Lemma~\ref{lem:B2} is satisfied and hence we
have
$$
\sum_{\lambda_r\in A_r'} \Bigl|\sum_{\lambda_1\in A_1}\ldots
\sum_{\lambda_{r-1}\in A_{r-1}} e^{2\pi i b \lambda_1\ldots \lambda_{r-1} \lambda_r/m_1}\Bigr|<m^{-\tau}|A_1|\ldots |A_r|
$$
for some absolute constant $\tau>0$ (see the discussion followed to
Lemma~\ref{lem:B2}). Inserting this into~\eqref{eqn:S1^4<} and using
estimates $k_i\ll (\log\log\log m)^{3}$ and $|A_i|\le M_i^{k_i},$ we
get
$$
|S_1|^4<m^{-\tau/5} M_1^{k_1}M_2^{k_2}\ldots M_r^{k_r}.
$$
Thus, from~\eqref{eqn:Case 1 S via S1} it follows that
$$
\frac{|S|}{M_1M_2\ldots M_r}<m^{-c_1(\log\log\log m)^{-3r}}
$$
and from~\eqref{eqn:W^2 Kloost} we get
$$
W<2 N^{1-0.5\beta}.
$$
Inserting this into~\eqref{eqn:sum x in G Kloost} and
using~\eqref{eqn:sum short x<N}, we conclude the proof.

\smallskip

{\it Case 2}. Let now $M_r < N^{\alpha^3}$. In this case we fix all
the factors except $p_1,p_2,p_r$. We apply
Corollary~\ref{cor:doubleKloostWithPrimes} or
Corollary~\ref{cor:doubleKloostWithAlmostPrimes}. We either choose
for the first factor $p_1$ and the second factor $p_2$ or for the
first factor $p_1p_r$ and the second factor $p_2$.
%Thus in
%Corollaries~\ref{cor:doubleKloostWithPrimes}
%and~\ref{cor:doubleKloostWithAlmostPrimes}
%the number $N_1$ is either $M_1$ or $M_1M_r$ and $N_2=M_2$.
Because $M_1 >
N^{\alpha}$ and $M_r <N^{\alpha^3}$ we will get in one of the cases
the required saving. Let us give some details of this argument.

Define $k_1,k_2\in \Z_{+}$ such that
$$
M_1^{k_1-1}<m_1^{1/2}\le M_1^{k_1},\qquad M_2^{k_2-1}<m_1^{1/2}\le M_2^{k_2}.
$$
From the definition of $\alpha$ and $\beta$ we have
$$
k_1\le \frac{1}{c}\log\log m;\qquad k_2\le
\frac{1}{c\beta}\ll \frac{(\log m)^{1/2}}{\log\log m}.
$$
Let
$$
\delta=\frac{k_1\log M_r}{3\log M_1}.
$$
Note that $\delta\le \frac{1}{c}(\log\log m)^{-1}.$ We further consider three subcases:

\smallskip

{\it Case 2.1.} Let $M_1^{k_1-1+\delta}<m_1^{1/2}\le
M_1^{k_1-\delta}$. Then we apply
Corollary~\ref{cor:doubleKloostWithPrimes} and get
\begin{eqnarray*}
\begin{split}
\frac{|S|}{M_1M_2\ldots M_r}& < (\log m)^{1/2}\Bigl(\frac{M_1^{k_1-1}}{m_1^{1/2}}+\frac{m_1^{1/2}}{M_1^{k_1}}\Bigr)^{1/(2k_1k_2)}\\ & < 2(\log m)^{1/2}M_1^{-\delta/(2k_1k_2)}
=2(\log m)^{1/2}M_r^{-1/(6k_2)}.
\end{split}
\end{eqnarray*}
Using the upper bound for $k_2$ and the lower bound $M_r\ge
N^{\beta}$ it follows that
$$
\frac{|S|}{M_1M_2\ldots M_r}<2(\log m)^{1/2} e^{-0.01c^2\beta^2\log m}.
$$

\smallskip

{\it Case 2.2.} Let $M_1^{k_1-\delta}<m_1^{1/2}\le M_1^{k_1}$. We
apply Corollary~\ref{cor:doubleKloostWithAlmostPrimes} in the form
\begin{eqnarray*}
\frac{|S|}{M_1M_2\ldots M_r}<(\log m)\Bigl(\frac{(M_1M_r)^{k_1-1}}{m_1^{1/2}}+\frac{m_1^{1/2}}{(M_1M_r)^{k_1}}\Bigr)^{1/(2k_1k_2)}\\
<(\log m)\Bigl(\frac{M_r^{k_1-1}}{M_1^{1-\delta}}+\frac{1}{M_r^{k_1}}\Bigr)^{1/(2k_1k_2)}.
\end{eqnarray*}

\smallskip

{\it Case 2.3.} Let now $M_1^{k_1-1}<m_1^{1/2}\le
M_1^{k_1-1+\delta}$.  Then $k_1\ge 2$ and we apply
Corollary~\ref{cor:doubleKloostWithAlmostPrimes} with $k_1$ replaced
by $k_1-1$  in the form
\begin{eqnarray*}
\frac{|S|}{M_1M_2\ldots M_r}<(\log m)^{1/2}\Bigl(\frac{(M_1M_r)^{k_1-2}}{m_1^{1/2}}+\frac{m_1^{1/2}}{(M_1M_r)^{k_1-1}}\Bigr)^{1/(2k_1k_2)}\\
<(\log m)^{1/2}\Bigl(\frac{M_r^{k_1-2}}{M_1}+\frac{M_1^{\delta}}{M_r^{k_1-1}}\Bigr)^{1/(2k_1k_2)}.
\end{eqnarray*}

\smallskip

In all three subcases we get the bound
$$
\frac{|S|}{M_1M_2\ldots M_r}<2(\log m)e^{-c'\beta^2\log m}
$$
for some constant $c'>0$. Thus, we eventually arrive at the bound
$$
W<N e^{-c''\beta^2\log m}\log m
$$
for some constant $c'>0$. Inserting this into~\eqref{eqn:sum x in G
Kloost} and using~\eqref{eqn:sum short x<N}, we conclude that
$$
\Bigl|\sum_{x<N}e_m(ax^*)\Bigr|\ll \beta(\log\log m)^{r-1} N+ N e^{-c'''\beta^2\log m}\log m\ll \frac{(\log\log m)^{r}}{(\log m)^{1/2}}N.
$$

\end{document}